\documentclass{amsart}
\usepackage{amsfonts}
\usepackage{amssymb}
\usepackage{amsmath}
\usepackage{amsthm}
\usepackage{hyperref}

\newcommand\Y{\mathbb Y}
\newcommand\Z{\mathbb Z}
\newcommand{\CC}{\mathbb C}
\newcommand\R{\mathbb R}

\newcommand\la{\lambda}

\newcommand\X{\frak X}

\newcommand\E{\bold{E}}
\newcommand\F{\bold{F}}
\newcommand\HH{\bold H}

\newcommand\wt{\widetilde}

\DeclareMathOperator*\const{{const}}

\DeclareMathOperator*\conf{{Conf}}

\DeclareMathOperator*\prob{{Prob}}

\DeclareMathOperator*\Projj{{Proj}}

\newcommand\unX{\underline X\,}
\newcommand\unK{\underline K\,}

\newcommand\g{{\textrm{gamma}}}
\newcommand\m{{\textrm{Meixner}}}
\newcommand\kr{{\textrm{Krawtchouk}}}
\newcommand\airy{{\textrm{Airy}}}

\newcommand\Bessel{{\textrm{Bessel}}}
\newcommand\sine{{\textrm{sine}}}

\newcommand\zz{{z,z',\xi}}

\newtheorem{Theorem}{Theorem}[section]
\newtheorem{Thesis}[Theorem]{Thesis}
\newtheorem{Corollary}[Theorem]{Corollary}

\newtheorem{Proposition}[Theorem]{Proposition}

\theoremstyle{definition}

\newtheorem{Remark}[Theorem]{Remark}

\numberwithin{equation}{section}

\begin{document}

\title{Difference operators and determinantal point processes}

\thanks{Supported by the RFBR grants
07-01-91209 and 08-01-00110, and by the research project SFB 701, University of
Bielefeld.}

\author{Grigori Olshanski}

\address{Institute for Information Transmission Problems. Bolshoy Karetny
19, Moscow 127994, GSP--4, Russia}

\email{olsh2007@gmail.com}

\dedicatory{Dedicated to I.~M.~Gelfand on the occasion of his 95th birthday}

\maketitle

\section{Introduction}

A point process is an ensemble of random locally finite point configurations in
a space $\X$. A convenient tool for dealing with point processes is provided by
the correlation functions (they are similar to the moments of random variables
\cite{Le}). The correlation functions are indexed by the natural numbers
$n=1,2,\dots$, and the $n$th function depends on $n$ variables ranging over
$\X$. A point process is said to be determinantal \cite{BO-CMP}, \cite{So} if
its correlation functions of all orders can be represented as the principal
minors of a correlation kernel $K(x,y)$, a function on $\X\times\X$. Then the
whole information about the point process is contained in this single function
in two variables.

Concrete examples of determinantal point processes arise in random matrix
theory, in some problems of representation theory, in combinatorial models of
mathematical physics, and in other domains. A remarkable fact is that the
correlation kernels of determinantal point processes from various sources
reveal a similar structure.

We will be interested in limit transitions for point processes depending on a
parameter. For instance, a typical problem in random matrix theory consists in
the study of limit properties of the spectra of random matrices of order $N$ as
the parameter $N$ goes to infinity. In the models of point processes arising in
representation theory, the limit transitions may be related to the
approximation of an infinite--dimensional group by finite--dimensional ones
\cite{Ol-PartI}, \cite{BO-AnnMath}.

The language of correlation functions is well adapted to studying limit
transitions: the convergence of processes is controlled by the convergence of
correlation functions (just as the convergence of random variables is
controlled by the convergence of their moments). For determinantal processes,
the situation is simplified, because the convergence of the correlation
functions is ensured by the convergence of the correlation kernels. In the
present note we will study just limit transitions for kernels.

Certain correlation kernels $K(x,y)$, which arise in many concrete examples
from various domains, share the following common properties:

\begin{itemize}

\item The space $\X$, on which the kernel is defined, is a subset of $\R$.

\item  The kernel can be written in the so--called integrable form,
\footnote{\,In the sense of Its, Izergin, Korepin, and Slavnov, see \cite{De}.}
which in the simplest variant looks as follows:

\begin{equation}\label{0.A}
K(x,y)=\dfrac{A(x)B(y)-B(x)A(y)}{x-y}\,,
\end{equation}
where $A$ and $B$ are certain functions on $\X$.

\end{itemize}

In such a situation, the convergence of kernels depending on a parameter is
often extracted from the asymptotics of the functions $A$ and $B$ with respect
to the parameter. This obvious way, however, may require hard computations if
$A$ and $B$ are sufficiently complex special functions. \footnote{For instance,
a solution of the problem of harmonic analysis on the infinite--dimensional
unitary group \cite{BO-AnnMath} required computing the asymptotics of of Hahn
type orthogonal polynomials, which are expressed through the hypergeometric
series ${}_3F_2$.} In the present note, we discuss another approach, which
exploits one more property of the kernels, which also holds in many concrete
examples:

\begin{itemize}

\item The operator $K$ in the Hilbert space $L^2(\X)$ (with respect to a
natural measure), given by the kernel $K(x,y)$, is a projection operator.
Moreover, $K$ can be realized as a spectral projection $\mathcal{P}(\Delta)$
for a certain selfadjoint operator $D$ acting in $L^2(\X)$. Here $\Delta$ is a
certain part of the spectrum of the operator $D$, and this operator is
determined by a difference or differential operator on $\X$.

\end{itemize}

As soon as there is a natural link between $K$ and $D$, the following simple
idea arises: to deduce the convergence of the correlation kernels from that of
the corresponding selfadjoint operators $D$. This idea was stated in
\cite{BO-JAlg} and then employed in  \cite{Go} and \cite{BG}. The aim of the
present note is to demonstrate the efficiency of such an approach on a number
of other examples. I consider only rather simple examples where the point
processes live on the one--dimensional lattice and not on a continuous space,
and the limit transition does not require a scaling. I think, however, that the
method can be useful in more complex situations, too. Some indication on this
is contained in \S3.4 where it is explained how one can guess, in a very simple
way, the scaling leading to the Airy kernel.

The structure of the work is as follows.  In \S1 we recall general definitions
related to point processes (for more detail, see \cite{Le}, \cite{So}). In \S2
we introduce the main model: a 3--parameter family of probability measures on
the Young diagrams  \footnote{These measures, called the z--measures, were
introduced in \cite{BO-CMP} for solving the problem of harmonic analysis on the
infinite symmetric group. See also the survey paper \cite{Ol-LN}.}). There it
is also explained how to pass from measures on the Young diagrams to point
processes on the lattice and where our main object, the ``hypergeometric''
difference operator on the lattice, comes from.  In \S3 we discuss limit
transitions related to the degeneration of the z--measures to the Plancherel
measure. In \S4 we study another limit regime that leads to the so--called
Gamma kernel. This interesting kernel (it arose in  \cite{BO-Gamma}) describes
the asymptotics of the fluctuations of the boundary of a random Young diagram
near the point of intersection with the diagonal  (in another language, limit
properties of the lowest Frobenius coordinates).

I am grateful to Alexei Borodin for discussions. I am also grateful to Rais
Ismagilov and Yuri Neretin for valuable critical remarks.

\section{Definitions}

Let $\X$ be a topological space. By a {\it configuration\/} in $\X$ we mean a
subset $X\subset\X$ without accumulation points. \footnote{More precisely, what
we termed a configuration should be called a simple locally finite
configuration. Here ``simple'' means that we exclude multiple points which are
allowed in a more general definition.} The space of all configurations in $\X$
is denoted as $\conf(\X)$. If $\X$ is discrete then $\conf(\X)$ is simply the
set $2^{\X}$ of all subsets of $\X$. We will be dealing with probability Borel
measures on $\conf(\X)$. Given such a measure $P$, one may speak about an
ensemble of {\it random\/} configurations or a {\it random point process\/}.

In what follows the space $\X$ will be discrete, so that all the next
definitions will be given for this simple case. About the general case, see,
e.g., \cite{Le}.

We assign to $P$ a sequence of functions on $\X$, $\X\times\X$,
$\X\times\X\times\X$, \dots, called the {\it correlation functions\/}. Here, by
definition, the $n$th correlation function $\rho^{(n)}$ is obtained in the
following way: its value $\rho^{(n)}(x_1,\dots,x_n)$ at an $n$--tuple of points
equals $0$ if among these points there are repetitions; otherwise the value is
equal to the probability of the event that the random configuration  $X$
contains all the points $x_1,\dots,x_n$.

The initial measure $P$ is determined by its correlation functions uniquely.

A point process $P$ is said to be {\it determinantal\/} \cite{BO-CMP},
\cite{So}, if there exists a function $K(x,y)$ on $\X\times\X$ such that for
any $n=1,2,\dots$
$$
\rho^{(n)}(x_1,\dots,x_n)=\det[K(x_i,x_j)]_{i,j=1}^n.
$$

Let us call $K(x,y)$ the {\it correlation kernel\/} of the process $P$. The
operator in the coordinate Hilbert space $\ell^2(\X)$ with the matrix
$[K(x,y)]$ will be called the  {\it correlation operator\/} and denoted by $K$.

The whole information on the determinantal point process $P$ is contained in
its correlation kernel (equivalently, correlation operator).

A simple but important example of determinantal point processes is afforded by
the {\it orthogonal polynomial ensembles\/}. Let us give the definition for the
discrete case we are interested in (for more details, see the survey paper
\cite{Ko}).

Let $\frak X\subset\R$ be a discrete subset and $W(x)>0$ be a weight function
defined on $\frak X$. Fix $N=1,2,\dots$. By definition, the ensemble in
question is formed by the $N$--point configurations
$X=\{x_1,\dots,x_N\}\subset\frak X$, whose probabilities are given by the
formula

\begin{equation}\label{1.A}
\prob(X)=\const\cdot\prod_{i=1}^N W(x_i)\cdot\prod_{1\le i<j\le N}(x_i-x_j)^2,
\end{equation}
where ``$\const$'' is a suitable normalizing factor.

It is well known (see, e.~g., \cite{Ko}) that such an ensemble is a
determinantal process and its correlation kernel  $K(x,y)$ is expressed through
the orthogonal polynomials $p_0=1,p_1,p_2,\dots$ with the weight function
$W(x)$, as follows:
$$
K(x,y)=\sum_{i=0}^{N-1}\wt p_i(x)\wt p_i(y), \qquad x,y\in\frak X,
$$
where
$$
\wt p_i(x)=\sqrt{W(x)}\,\frac{p_i(x)}{\Vert p_i\Vert_W}\,, \qquad x\in\frak X,
\quad i=0,1,\dots,
$$
and $\Vert \,\cdot\,\Vert_W$ stands for the norm in the weighted Hilbert space
$\ell^2(\frak X,W)$. Observe that the functions $\wt p_i(x)$ form an
orthonormal system of vectors in $\ell^2(\frak X)$. Thus, the kernel $K(x,y)$
is the matrix of a finite--dimensional projection in $\ell^2(\frak X)$:
specifically, of the projection onto the subspace of polynomials of degree at
most $N-1$, multiplied by $\sqrt{W(x)}$.

Since $K(x,y)$ is equal to the product of  $\sqrt{W(x)W(y)}$ with the $N$th
Christoffel--Darboux kernel for our system of orthogonal polynomials, the
kernel  $K(x,y)$ can be written in the form
$$
K(x,y)=\const\cdot\frac{\wt p_N(x)\wt p_{N-1}(y)-\wt p_{N-1}(x)\wt p_N(y)}{x-y}
$$
and hence can be represented in the integrable form (\ref{0.A}).

\section{Z--measures and discrete orthogonal polynomial ensembles}

\subsection{Z--measures \cite{BO-CMP}, \cite{BO-MMJ}}
Denote by $\Y$ the set of all partitions $\la=(\la_1,\la_2,\dots)$, which we
identify with the corresponding Young diagrams. The {\it z--measure\/} with
parameters
$$
z\in\CC, \qquad z'\in\CC, \qquad\xi\in\CC\setminus[1,+\infty)
$$
is the (complex) measure $M_\zz$ on the set $\Y$ assigning to a diagram
$\la\in\Y$ the weight

\begin{equation}\label{2.F}
M_\zz(\la)=(1-\xi)^{zz'}\,\xi^{|\la|}\,(z)_\la(z')_\la\,
\left(\frac{\dim\la}{|\la|!}\right)^2\,.
\end{equation}
Here $|\la|:=\sum_i\la_i$ (equivalently, $|\la|$ is the number of boxes in the
diagram  $\la$);
\begin{equation}\label{2.J}
(z)_\la=\prod_{i=1}^{\ell(\la)}(z-i+1)_{\la_i}=\prod_{(i,j)\in\la}(z+j-i),
\end{equation}
where $(a)_m=a(a+1)\dots(a+m-1)$ is the Pochhammer symbol; $\ell(\la)$ is the
number of nonzero coordinates $\la_i$; the second product in (\ref{2.J}) is
taken over all boxes  $(i,j)$ of the diagram $\la$, where $i$ and $j$ denote
the numbers of the row and the column containing the box; finally, $\dim\la$
denotes the number of standard tableaux of shape $\la$ (equivalently, the
dimension of the irreducible representation of the symmetric group of degree
$|\la|$ indexed by the diagram $\la$).

The following summation formula holds:
$$ \sum_{\la\in\Y}M_\zz(\la)=1.
$$
Here we assume that either $|\xi|<1$ (then the series absolutely converges) or
the parameters $z,z'$ are integers of opposite sign (then the series
terminates).

Note that the z--measures are a particular case of the Schur measures, see
\cite{Ok-Infw}).

Note also two properties of the z--measures:

\begin{equation}\label{2.I}
M_\zz(\la)=M_{z',z,\xi}(\la) \quad \textrm{и} \quad
M_\zz(\la')=M_{-z,-z',\xi}(\la),
\end{equation}
where $\la'$ denotes the transposed diagram.

In what follows we consider only the z--measures with real nonnegative weights:
$M_\zz(\la)\ge0$ for all $\la\in\Y$. These are probability measures on $\Y$.
They fall into the following 4 series.

\medskip
$\bullet$ {\it Principal series\/}: The parameters $z$ and $z'$ are complex
conjugate and nonreal while the parameter $\xi$ is real and $0<\xi<1$.
\medskip

$\bullet$ {\it Complementary series\/}: The both parameters $z$ and $z'$ are
real and contained in one and the same open interval of the form $(m,m+1)$
where $m\in\Z$, while the parameter $\xi$ is the same as above ($0<\xi<1$).
\medskip

$\bullet$ {\it First degenerate series\/}: One of the parameters $z,z'$ (say,
$z$) is a nonzero integer, while the other parameter (then it will be $z'$) is
a real number of the same sign such that $|z'|>|z|-1$. Again, $0<\xi<1$.
\medskip

$\bullet$ {\it Second degenerate series\/}: The parameters $z$ and $z'$ are
integers of opposite sign. The parameter $\xi$ is now subject to another
condition: $\xi<0$.
\medskip

If a  z--measure belongs to the principal or complementary series then
$(z+k)(z'+k)>0$ for all integers $k\in\Z$, which implies that the weights of
all diagrams are strictly positive. For the z--measures of the degenerate
series some of the weights vanish but all nonzero weights are strictly
positive.

\subsection{More details on degenerate z--measures. The Meixner and Krawtchouk ensembles}

Consider the first degenerate series. By virtue of the symmetry relation
(\ref{2.I}) it suffices to take $z=N$ and $z'=N+c-1$, where $N=1,2,\dots$ and
$c>0$. Then the weight of a diagram $\la$ vanishes precisely when
$\ell(\la)>N$. Thus, the support of the z--measure is the set of the diagrams
contained in the horizontal strip of width $N$. Such diagrams $\la$ are in a
one--to--one correspondence with the $N$--point configurations $L$ on
$\Z_+:=\{0,1,2,\dots\}$:

\begin{equation}\label{2.A}
\la=(\la_1,\dots,\la_N,0,0,\dots)\; \longleftrightarrow \; L=(l_1,\dots,l_N),
\end{equation}
where
$$
l_i=\la_i+N-i, \quad 1\le i\le N.
$$
It is readily checked that in terms of the correspondence (\ref{2.A}), the
weight $M_{N,N+c-1,\xi}(\la)$ can be written in the form (\ref{1.A}), where we
have to set $x_i=l_i$ and take as  $W$ the weight function for the {\it Meixner
polynomials\/} \cite{KS} on the set $\frak X=\Z_+$:
$$
W^\m(l)=\frac{(c)_l\xi^{ l}}{l!}\,, \qquad   l\in\Z_+\,.
$$
Thus, the correspondence $\la\to L$ converts the z--measure of the first
degenerate series with parameters $z=N$, $z'=N+c-1$, $\xi$ into the $N$--point
Meixner orthogonal polynomial ensemble with parameters $c$ and $\xi$.

A similar fact holds for the second degenerate series. Let $z=N$ and $z'=-N'$
where $N$ and $N'$ are two positive integers, and assume $\xi<0$. Then the
weight of a diagram $\la$ does not vanish if and only if $\ell(\la)\le N$ and
$\ell(\la')\le N'$, that is, $\la$ has to be contained in the rectangular
$\square_{N,N'}$ with $N$ rows and $N'$ columns. The support of such a measure
is finite.

The same correspondence  $\la\leftrightarrow L$ as in (\ref{2.A}) gives a
bijection between the diagrams $\la\subseteq\square_{N,N'}$ and the $N$--point
configurations $L$ on the finite set $\{0,1,\dots,\wt N\}$, where $\wt
N=N+N'-1$. Then the z--measure of the second degenerate series turns into the
$N$--point orthogonal polynomial ensemble on $\frak X=\{0,1,\dots,\wt N\}$ that
is determined by the weight function
$$
W^\kr(l)= \binom{\wt N}{x}p^l(1-p)^{\wt N-l}, \qquad  l=0,\dots,\wt N,
$$
for the {\it Krawtchouk orthogonal polynomials\/} \cite{KS} with the parameters
$$
p:=\frac\xi{\xi-1}\,, \qquad \wt N:=N+N'-1.
$$
Note that the condition $\xi<0$ guarantees that $0<p<1$.

Thus, the Meixner and Krawtchouk ensembles can be extracted from the
z--measures as a rather particular case. As we will see, it is also helpful to
use the inverse transition:

\begin{Thesis}\label{Thesis}
The z--measures can be obtained from the Meixner or Krawtchouk ensembles via
analytic continuation (interpolation) with respect to the parameters.
\end{Thesis}

Let us explain this statement (for more detail, see \cite{BO-MMJ} and
\cite{BO-Markov}). Note that the support of a degenerate measure enlarges as
the parameter $N$ grows (or the two parameters $N$ and $N'$ grow), so that any
given diagram $\la\in\Y$ can be covered. Further, as seen from (\ref{2.F}), for
any fixed diagram $\la$, its weight $M_\zz(\la)$ is an analytic function in
$\xi$ whose Taylor coefficients at $\xi=0$ are given by polynomial functions in
$z$ and $z'$. This allows one to extrapolate results about the degenerate
series to the general case, because such functions in  $z$, $z'$, and $\xi$ are
uniquely determined by their restriction to the subset of the values of the
parameters corresponding to the degenerate series.

As we will see, Thesis 2.1 explains the origin of the difference operator
(\ref{2.C}) introduced below.

\subsection{The hypergeometric deifference operator}

Consider now the principal and complementary series. We cannot use the
correspondence (\ref{2.A}) anymore. Instead of it we will introduce another
correspondence, which will take {\it arbitrary\/} partitions $\la\in\Y$ to {\it
semi--infinite\/} point configurations $\unX$ on the lattice $\Z'=\Z+\frac12$
of half--integers:

\begin{equation}\label{2.B}
\la\;\longleftrightarrow\;\unX=\{\la_i-i+\tfrac12\mid i=1,2,\dots\}.
\end{equation}

(We call the configuration  $\unX$ semi--infinite because it contains all the
lattice points that are sufficiently far to the left of zero and does not
contain points sufficiently far to the right of zero.)

There is a simple link between (\ref{2.A}) and  (\ref{2.B}). Given a diagram
$\la\in\Y$, let $N$ be so large that $N\ge\ell(\la)$, so that the $N$--point
configuration $L$ determined by (\ref{2.A}) exists. Then $\unX$ is obtained
from  $L$ by shifting it to the left by $N-\frac12$ and next adding the
left--infinite ``tail'' $\{\tfrac12-i\mid i=N+1,N+2,\dots\}$.

Let us assume the parameters $(z,z')$ lie in the principal or complementary
series. Consider the following difference operator on the lattice $\Z'$:

\begin{multline}\label{2.C}
\mathcal{D}_\zz f(x)=\sqrt{\xi(z+x+\tfrac12)(z'+x+\tfrac12)}\,
f(x+1)-(x+\xi(z+z'+x))f(x)\\+\sqrt{\xi(z+x-\tfrac12)(z'+x-\tfrac12)}\, f(x-1).
\end{multline}

This difference operator was introduced in the papers \cite{BO-Markov} and
\cite{BO-MMJ} as the result of analytic continuation, in the spirit of Thesis
\ref{Thesis}, of the Meixner difference operator. In more details, the analytic
continuation procedure is as follows:

\medskip

1. Consider the Meixner difference operator on $\Z_+$:
$$
\mathcal{D}^\m_{c,\xi}f(x)=\xi(x+c)f(x+1)-(x+\xi(x+c))f(x)+xf(x-1).
$$
It multiplies the $n$th Meixner polynomial by $-(1-\xi)n$, see \cite{KS},
formula (1.9.5). Add to $\mathcal{D}^\m_{c,\xi}$ the constant term
$(1-\xi)(N-\frac12)$ and observe that the first $N$ Meixner polynomials are
just those eigenvectors of the resulting operator

\begin{equation}\label{2.D}
\mathcal{D}^\m_{c,\xi}+(1-\xi)(N-\tfrac12)
\end{equation}
that correspond to nonnegative eigenvalues.
\medskip

2. Pass from the weight $\ell^2$--space $\ell^2(\Z_+,W^\m)$ to the ordinary
space $\ell^2(\Z_+)$. This means that (\ref{2.D}) should be replaced by the
composition of operators

\begin{equation}\label{2.E}
 (W^\m)^{\frac12}\circ\left(\mathcal{D}^\m_{c,\xi}+(1-\xi)(N-\tfrac12)\right)\circ(W^\m)^{-\frac12}.
\end{equation}
\medskip

3. Make in (\ref{2.E}) the change of the variable  $l\to x:=l-(N-\frac12)$,
which means that the support $\Z_+$ of the weight function is transformed to
the subset $\{-(N-\frac12)+m\mid m=0,1,2,\dots\}\subset\Z'$. Note that this
subset grows together with  $N$ and in the limit it exhausts the whole lattice
$\Z'$.
\medskip

4. Formally substitute $N=z$ and $c=z'-z+1$ into the coefficients of the
resulting operator.
\medskip

Then we get $\mathcal{D}_\zz$. Likewise, $\mathcal{D}_\zz$ can be obtained from
the difference operator connected to the Krawtchouk polynomials.

Denote by $\{e_x\}$, $x\in\Z'$, the natural orthonormal basis in the Hilbert
space  $\ell^2(\Z')$ and let $\ell^2_0(\Z')\subset\ell^2(\Z')$ stand for the
algebraic subspace formed by linear combinations of the basis elements. We will
consider  $\mathcal{D}_\zz$ (see (\ref{2.C})) as a symmetric operator in
$\ell^2(\Z')$ with domain $\ell^2_0(\Z')$.

\begin{Proposition}\label{2.3}
The operator $\mathcal{D}_\zz$ defined in this way is essentially selfadjoint.
\end{Proposition}

\begin{proof}
The standard way is to check that the eigenfunctions of the difference operator
 (\ref{2.C}) with nonreal eigenvalues do not belong to
$\ell^2(\Z')$. This can be done by computing explicitly the eigenfunctions
(they are expressed through the Gauss hypergeometric function, see
\cite{BO-MMJ}). The computation can be simplified if one replaces  (\ref{2.C})
with the difference operator
$$
\sqrt\xi\,|x+\tfrac12|\, f(x+1)-(x+\xi(z+z'+x))f(x)+\sqrt\xi\,|x-\tfrac12|\,
f(x-1),
$$
which differs from (\ref{2.C}) by a bounded operator (adding a bounded
selfadjoint operator preserves essential selfadjointness).

I will sketch now another argument. It is perhaps less elementary but it
requires no computations and builds a bridge to representations of the group
$SL(2)$, which turns out to be useful in other questions related to the
z--measures.

Consider the complex Lie algebra $\frak{g}_\CC=\frak{sl}(2,\CC)$ with the basis
$$
\E=\bmatrix 0 & 1\\0 & 0 \endbmatrix, \quad \F=\bmatrix 0 & 0\\1 & 0
\endbmatrix, \quad \HH=\bmatrix 1 & 0\\0 & -1 \endbmatrix
$$
and the commutation relations
$$
[\HH,\E]=2\E, \quad [\HH,\F]=-2\F, \quad [\E,\F]=\HH.
$$
The real linear span of the elements $\E+\F$, $i(\E-\F)$, and $i\HH$ is a real
form $\frak g\subset\frak{g}_\CC$: this is the Lie algebra $\frak{su}(1,1)$.

Define a representation $S_{z,z'}$ of the Lie algebra $\frak{g}_\CC$ in the
algebraic space $\ell^2_0(\Z')\subset\ell^2(\Z')$ by the following formulas for
the basis elements

\begin{equation}\label{2.H}
\gathered S_{z,z'}(\E)\,e_x=\sqrt{(z+x+\tfrac12)(z'+x+\tfrac12)}\;e_{x+1}\\
S_{z,z'}(\F)\,e_x=-\;\sqrt{(z+x-\tfrac12)(z'+x-\tfrac12)}\;e_{x-1}\\
S_{z,z'}(\HH)\,e_x=(z+z'+2x)e_x
\endgathered
\end{equation}
It is readily checked that the above operators satisfy the commutation
relations and hence do determine a representation.

The representation $S_{z,z'}$ equips $\ell^2_0(\Z')$ with the structure of an
irreducible Harish--Chandra module over $\frak g$ such that this algebra acts
by skew symmetric operators. It follows that $S_{z,z'}$ gives rise to an
irreducible unitary representation of the group $SU(1,1)^\sim$, the simply
connected covering group of the matrix group $SU(1,1)=SL(2,\R)$ (see \cite{Di},
p.~14, where a proof is given for a concrete example but the argument holds in
the general case). Moreover, according to a well--known theorem by
Harish--Chandra, all vectors from  $\ell^2_0(\Z')$ are analytic vectors of this
representation.

On the other hand, as seen from (\ref{2.C}) and (\ref{2.H}),

\begin{equation}\label{2.G}
\mathcal{D}_\zz=S_{z,z'}\left(\sqrt\xi \,\E-\sqrt\xi
\,\F-\tfrac{1+\xi}2\,\HH\right)+\tfrac{1-\xi}2(z+z')\bold1,
\end{equation}
where $\bold1$ denotes the identity operator. Therefore, all vectors from
$\ell^2_0(\Z')$ are analytic vectors of the operator $\mathcal{D}_\zz$, which
implies (\cite{RS}, теорема X.39) that it is essentially selfadjoint.
\end{proof}

\begin{Remark}\label{2.4} The use of the representation $S_{z,z'}$ of the Lie
algebra $\frak{sl}(2,\CC)$ is suggested by Okounkov's work \cite{Ok-SL2}. Note
that the corresponding unitary representation of the group $SU(1,1)^\sim$
belongs to the principal or complementary series, in exact agreement with the
terminology for the z--measures (which was introduced prior to \cite{Ok-SL2},
relying on a purely formal analogy).
\end{Remark}

Let us agree to denote by $D_\zz$ the selfadjoint operator in $\ell^2(\Z')$
that is the closure of the operator $\mathcal{D}_\zz$. It follows from
Proposition \ref{2.3} that the domain of the operator $D_\zz$ consists of all
those functions $f\in \ell^2(\Z')$ that remain in $\ell^2(\Z')$ after
application of the difference operator (\ref{2.C}).

\begin{Proposition}\label{2.5}
The selfadjoint operator $D_\zz$ has simple, pure discrete spectrum filling the
subset $(1-\xi)\Z'\subset\R$.
\end{Proposition}

\begin{proof}
The matrix $\sqrt\xi \,\E-\sqrt\xi \,\F-\tfrac{1+\xi}2\,\HH$ that enters
(\ref{2.G}) is conjugated with the matrix $-\frac{1-\xi}2\HH$ by means of an
element of $SU(1,1)$. Therefore, $D_\zz$ is unitarily equivalent to the closure
of the operator $S_{z,z'}(-\frac{1-\xi}2\HH)+\tfrac{1-\xi}2(z+z')\bold1$. This
operator is the diagonal operator in the basis $\{e_x\}$ acting as
multiplication of $e_x$ by $-(1-\xi)x$. This makes the claim evident.
\end{proof}

Thus, for any triple of parameters $(\zz)$ from the principal or complementary
series there exists an orthonormal basis $\{\psi_{a;\zz}\}$ in $\ell^2(\Z')$,
indexed by points $a\in\Z'$ and such that
$$
\mathcal{D}_\zz \psi_{a;\zz}=(1-\xi)a\,\psi_{a;\zz}, \qquad a\in\Z'.
$$
As shown in \cite{BO-MMJ}, the eigenvectors $\psi_{a;\zz}$ can be written down
explicitly: for them there exists a convenient contour integral representation
and also an expression through the Gauss hypergeometric function.

In view of the connection of the functions $\psi_{a;\zz}$ with the
hypergeometric function I will call  $\mathcal{D}_\zz$ the {\it hypergeometric
difference operator\/}.

If $A$ is a selfadjoint operator for which $0$ is not a point of discrete
spectrum, then we will denote by $\Projj_+(A)$ the spectral projection
corresponding to the positive part of the spectrum of $A$. In particular, as
the spectrum of $D_\zz$ does not contain $0$, we may form the projection
$\Projj_+(D_\zz)$.

The image of the $z$--measure $M_\zz$ under the correspondence $\la\mapsto\unX$
introduced in (\ref{2.B}) determines a point process on $\Z'$, which we denote
as $P_\zz$.

\begin{Theorem}\label{2.6} $P_\zz$ is a determinantal process. Its correlation
kernel, denoted as $\unK_{z,z',\xi}(x,y)$, is the matrix of the spectral
projection $\Projj_+(D_\zz)$. That is,
$$
\unK_\zz(x,y)=\sum_{a\in\Z'_+}\psi_{a;\zz}(x)\psi_{a;\zz}(y), \qquad x,y\in\Z'.
$$
\end{Theorem}

(We do not put the bar over $\psi_{a;\zz}(y)$ because the eigenfunctions are
real--valued, see formula (2.1) in \cite{BO-MMJ}.)

In the above formulation, the result is contained in \cite{BO-MMJ}. There are
also references to previous works. We call $K_{z,z',\xi}(x,y)$ the {\it
discrete hypergeometric kernel\/.} It can be written in the integrable form
(\ref{0.A}), see \cite{BO-MMJ}, Prop. 3.10.

\section{Plancherel measure, dicrete Bessel kernel, and discrete sine kernel}

\subsection{Poissonized Plancherel measure}

Consider the limit regime

\begin{equation}\label{3.A}
\xi\to0, \quad z\to\infty, \quad z'\to\infty, \quad \xi zz'\to\theta,
\end{equation}
where $\theta>0$ is a new parameter. Then the $z$--measures $M_\zz$ converge to
a probability measure $M_\theta$ on $\Y$:
$$
M_\theta(\la)=\lim_{zz\xi\to\theta}M_\zz(\la)
=e^{-\theta}\theta^{|\la|}\left(\frac{\dim\la}{|\la|!}\right)^2.
$$

The measure $M_\theta$ arises in the result of poissonization of the sequence
the Plancherel measures $M^{(n)}$, $n=1,2,\dots$, where
$$
M^{(n)}(\la)=\frac{(\dim\la)^2}{n!}, \qquad |\la|=n.
$$
For this reason $M_\theta$ is called the {\it poissonized Plancherel measure\/}
with parameter $\theta$.

For more details about the Plancherel measures $M^{(n)}$ see the papers
\cite{LS} and \cite{VK}; see also  \cite{Ke} and \cite{IO}. The poissonized
version $M_\theta$ was first considered in  \cite{BDJ}; further results were
obtained in \cite{BOO} and \cite{Jo}.

As in \S2, we use the correspondence $\la\mapsto\unX$ to pass from the measure
$M_\theta$ on $\Y$ to a point process on  $\Z'$; let us denote this point
process as $P_\theta$. It is determinantal; its correlation kernel is described
in \S3.2 below.

\subsection{{}From discrete hypergeometric kernel to discrete Bessel kernel}

In the regime (\ref{3.A}), the formal limit of the difference operator
$\mathcal{D}_\zz$ is the difference operator $\mathcal{D}^\Bessel_\theta$ on
the lattice $\Z'$, acting according to formula
$$
\mathcal{D}^\Bessel_\theta f(x)=\sqrt\theta f(x+1) -xf(x)+\sqrt\theta f(x-1).
$$

\begin{Proposition}\label{3.1}
Consider $\mathcal{D}^\Bessel_\theta$ as a symmetric operator in $\ell^2(\Z')$
with domain $\ell^2_0(\Z')$. Then this is an essentially selfadjoint operator.
\end{Proposition}

\begin{proof}
Indeed, $\mathcal{D}^\Bessel_\theta$ is the sum of the diagonal operator
$f\mapsto -xf$ and a bounded operator.
\end{proof}

Alike the hypergeometric operator (\ref{2.C}), the operator
$\mathcal{D}^\Bessel_\theta$ is connected with a group representation, only the
group is different, namely, it is the universal covering group $\wt G$ for the
group $G$ of motions of the plane $\R^2$. Let us discuss this connection
shortly.

Let $\frak{g}$ be the Lie algebra of the group $G$ and $\frak{g}_\CC$ be the
complexification of $\frak{g}$. In $\frak{g}_\CC$, there is a basis
$\{\E,\F,\HH\}$ with the commutation relations
$$
[\HH,\E]=\E, \quad [\HH,\F]=-\F, \quad [\E,\F]=0.
$$
Consider the representation $S$ of the Lie algebra $\frak{g}_\CC$ in the dense
subspace $\ell^2_0(\Z')\subset\ell^2(\Z')$ defined on the basis elements by the
formulas
$$
S(\E)\,e_x=e_{x+1}, \quad S(\F)\,e_x=-e_{x-1}, \quad S(\HH)e_x=xe_x,
$$
It is directly checked that all vectors from  $\ell^2_0(\Z')$ are analytic
vectors for $S$ and $\frak g$ acts by skew symmetric operators. Therefore, $S$
gives rise to a unitary representation of the group $\wt G$. We have
$$
\mathcal{D}^\Bessel_\theta=S\left(\sqrt\theta \,\E-\sqrt\theta \,\F-\HH\right),
$$
which is an analog of equality (\ref{2.G}).

Denote by  $D^\Bessel_\theta$ the closure of the symmetric operator
$\mathcal{D}^\Bessel_\theta$. According to Proposition  \ref{3.1}, $D_\theta$
is selfadjoint.

Set
$$
\psi_{a;\theta}(x):=J_{x+a}(2\sqrt\theta), \qquad x\in\Z',
$$
where $a\in\Z'$ is a parameter and  $J_m(\,\cdot\,)$ is the Bessel function of
index $m$.

\begin{Proposition}\label{3.2} For any fixed $\theta>0$, the functions
$\psi_{a;\theta}$ form an orthonormal basis in $\ell^2(\Z')$ and
$$
D^\Bessel_\theta\,\psi_{a;\theta}=a\,\psi_{a;\theta}, \qquad a\in\Z'.
$$
Thus, $D^\Bessel_\theta$ has pure discrete, simple spectrum filling the lattice
$\Z'\subset\R$.
\end{Proposition}

\begin{proof}
The element $\sqrt\theta \,\E-\sqrt\theta \,\F-\HH\in i\frak g$ is conjugated
with $-\HH$ under the adjoint action of the group $G$. It follows that the
selfadjoint operators $\overline{S(\sqrt\theta \,\E-\sqrt\theta \,\F-\HH)}$ and
$\overline{S(-\HH)}$ (where the bar means closure) are conjugated by means of a
unitary operator corresponding to an element of the group $\wt G$. Since
$D^\Bessel_\theta$ coincides with $\overline{S(\sqrt\theta \,\E-\sqrt\theta
\,\F-\HH)}$, its spectrum is the same as that of the diagonal operator
$\overline{S(-\HH)}$. Thus, we see that the spectrum is indeed the same as was
stated.

The fact that the functions  $\psi_{a;\theta}(x)$ are eigenfunctions of the
operator $\mathcal{D}^\Bessel_\theta$ follows from the well--known recurrence
relations for the Bessel functions (\cite{Er},7.2.8 (56)). Finally, the
relation $\Vert\psi_{a;\theta}\Vert^2=1$ can be proved in the same way as
Proposition 2.4 in \cite{BO-MMJ}, by using the contour integral representation
of the Bessel functions.
\end{proof}

\begin{Proposition}\label{3.3}
In the regime\/ {\rm(\ref{3.A})}, the operators $D_\zz$ converge to the
operator {\rm $D^\Bessel_\theta$} in the strong resolvent sense.
\end{Proposition}

About the notion of strong resolvent convergence see the textbook \cite{RS},
Section VIII.7.

\begin{proof}
Obviously, $D_\zz\to D^\Bessel_\theta$ on the subspace $\ell^2_0(\Z')$. By
virtue of Propositions  \ref{2.3} and \ref{3.1}, $\ell^2_0(\Z')$ is a common
essential domain for all selfadjoint operators under consideration. Then the
claim follows from a well--known general theorem (\cite{RS}, Theorem VIII.25).
\end{proof}

\begin{Corollary}\label{3.4}
The spectral projection operators $\Projj_+(D_\zz)$ strongly converge to the
spectral projection operator {\rm $\Projj_+(D^\Bessel_\theta)$}.
\end{Corollary}

\begin{proof}
This is a direct consequence of Proposition \ref{3.3} and another general
theorem (\cite{RS}, Theorem VIII.24 (b)). For applicability of that theorem it
is important that $0$ not be a point of discrete spectrum of the operators
$D_\zz$ and $D^\Bessel_\theta$, and this fact follows from the description of
their spectra, see Propositions  \ref{2.5} and \ref{3.2}.
\end{proof}

Corollary \ref{3.4} shows that $\Projj_+(D^\Bessel_\theta)$ serves as the
correlation operator for $P_\theta$.  The kernel of the projection
$\Projj_+(D^\Bessel_\theta)$ can be written in the form
$$
K^\Bessel_\theta(x,y)=\sum_{a\in\Z'_+}\psi_{a;\theta}(x)\psi_{a;\theta}(y)
=\sqrt\theta\,\frac{\psi_{-\frac12;\theta}(x)\psi_{\frac12;\theta}(y)
-\psi_{\frac12;\theta}(x)\psi_{-\frac12;\theta}(y)}{x-y}
$$
This kernel is called the {\it discrete Bessel kernel\/}; it was independently
derived in  \cite{Jo} and (in somewhat different form) in \cite{BOO}.

\subsection{{}From discrete Bessel kernel to discrete sine kernel}
Consider the following limit regime depending on the parameter $c\in(-1,1)$:

\begin{equation}\label{3.B}
\theta\to\infty, \qquad x\approx 2c\sqrt\theta+\wt x, \qquad  \wt x\in\Z.
\end{equation}
This means that as $\theta$ goes to $\infty$, we shift the lattice  $\Z'$ so
that to focus on a neighborhood of the point $2c\sqrt\theta$. The following
formal limit holds
$$
\left(\frac1{\sqrt\theta}\,\mathcal{D}^\Bessel_\theta\right)_{x\rightsquigarrow
\tilde x}\;\to\;\mathcal{D}^\sine_c
$$
with the following difference operator in the right--hand side:
$$
\mathcal{D}^\sine_c f(\wt x)=f(\wt x+1)-2cf(\wt x)+f(\wt x-1), \qquad  \wt
x\in\Z.
$$

\begin{Proposition}\label{3.5}
Consider {\rm $\mathcal{D}^\sine_c$} as a symmetric operator in $\ell^2(\Z)$
with domain  $\ell^2_0(\Z)$. Then this is an essentially selfadjoint operator.

Its closure, denoted as {\rm $D^\sine_c$}, has pure discrete spectrum of
multiplicity $2$, filling the interval $(-2-2c, 2-2c)$.
\end{Proposition}

\begin{proof}
The difference operator $\mathcal{D}^\sine_c$ is invariant under shifts of
$\Z$, which makes its study quite simple. Indeed, pass from the lattice $\Z$ to
the unit circle $|\zeta|=1$ by means of the Fourier transform. Then our
operator will become the operator of multiplication by the function
$$
\zeta\;\to\;2(\Re\zeta-c),
$$
which easily implies the claim about the spectrum.
\end{proof}

\begin{Proposition}\label{3.6} In the regime\/ {\rm (\ref{3.B})}, the operators
{\rm $D^\Bessel_\theta$} converge to the operator\/  {\rm $D^\sine_c$} in the
strong resolvent sense.
\end{Proposition}

The proof is the same as in Proposition \ref{3.3}. As in Corollary \ref{3.4},
we deduce from this the strong convergence $\Projj_+(D^\Bessel_\theta)\to
\Projj_+(D^\sine_c)$. The kernel of the projection $\Projj_+(D^\sine_c)$ is
readily computed: it has the form
$$
K^\sine_c(x,y)=\frac{\sin((\arccos c)(x-y))}{\pi(x-y)}, \qquad
 x,y\in\Z.
$$
and is called the {\it discrete sine kernel \/} \cite{BOO}.

\subsection{{}From discrete Bessel kernel to Airy kernel}
The Airy kernel first emerged in random matrix theory: under suitable
conditions on the matrix ensemble, this kernel describes the asymptotics of the
spectrum ``at the edge'', as the order of matrices goes to infinity. The Airy
kernel is defined on the whole real line and is expressed through the classical
Airy function $Ai(u)$ and its derivative:
$$
K^\airy(u,v)=\int_0^{+\infty}Ai(u+s)Ai(v+s)ds=\frac{Ai(u)Ai'(v)-Ai'(u)Ai(v)}{u-v},
\qquad u,v\in\R.
$$

As shown in \cite{BOO} and \cite{Jo}, the discrete Bessel kernel
$K^\Bessel_\theta(x,y)$ converges to the Airy kernel $K^\airy(u,v)$ if
$\theta\to+\infty$ and the variable $x\in\Z'$ is related to the variable
$u\in\R$ by the scaling $x=2\theta^{1/2}+\theta^{1/6}\cdot u$. Below we present
a heuristic derivation of this claim.

Let $g(u)$ be a smooth function on $\R$. Assign to it the function  $f(x)$ on
$\Z'$ by setting $f(x)=g(u)$, with the understanding that $x$ and $u$ are
related to each other as above. Then we get
$$
f(x\pm1)=g(u\pm\theta^{-1/6})\approx g(u)\pm\theta^{-1/6}\cdot
g'(u)+\tfrac12\,\theta^{-1/3}\cdot g''(u),
$$
so that
$$
\mathcal{D}^\Bessel_\theta
f(x)=\theta^{1/2}(f(x+1)-x\theta^{-1/2}f(x)+f(x-1))\approx\theta^{1/6}(g''(u)-ug(u)).
$$
This simple computation explains the origin of the factor $2$ and the exponents
$1/2$ and  $1/6$.

Next, the eigenvalue equation
$$
\mathcal{D}^\Bessel_\theta\psi=a\psi, \quad a\in\Z',
$$
turns, after the renormalization $a=\theta^{1/6}s$, into the equation
$$
\mathcal{D}^\airy\psi=s\psi, \quad s\in \R,
$$
where
$$
\mathcal{D}^\airy g(u)=g''(u)-ug(u)
$$
is the Airy operator. Its eigenfunctions have the form
$$
\psi_{s}(u)=Ai(u+s), \qquad s\in\R,
$$
and the Airy kernel corresponds to the spectral projection onto the positive
part of the spectrum.

It would be interesting to make this argument rigorous by using the notion of
Mosco convergence.

\section{{}From discrete hypergeometric kernel to Gamma kernel}

Let us fix the parameters  $(z,z')$ from the principal or complementary series
and let $\xi$ go to $1$. As shown in  \cite{BO-Gamma}, in this regime, there
exists a limit for the correlation kernels. Below we prove this fact in a new
way (Corollary \ref{4.4}).

The convergence of the kernels implies that the images of the measures  $M_\zz$
under the correspondence (\ref{2.B}) converge to a probability measure on
$\conf(\Z')$ generating a determinantal point process. Note that for the
initial  $z$--measures on $\Y$ the picture is different: they converge to $0$.
Of course, this cannot happen on the compact space $\conf(\Z')=2^{\Z'}$.

The first indication on the existence of a limit is the fact that the
hypergeometric difference operator (\ref{2.C}) still makes a sense under the
formal substitution $\xi=1$: then we get the difference operator

\begin{multline}\label{4.A}
\mathcal{D}^\g_{z,z'}f(x)=\sqrt{(z+x+\tfrac12)(z'+x+\tfrac12)}\,
f(x+1)-(z+z'+2x)f(x)\\+\sqrt{(z+x-\tfrac12)(z'+x-\tfrac12)}\, f(x-1).
\end{multline}
In the sequel we follow the scheme we have already worked out.

\begin{Proposition}\label{4.1}
Consider {\rm $\mathcal{D}^\g_{z,z'}$} as a symmetric operator in $\ell^2(\Z')$
with domain $\ell^2_0(\Z')$. Then this is an essentially selfadjoint operator.
\end{Proposition}

\begin{proof}
The both arguments given in the proof of Proposition \ref{2.3} hold without
changes.
\end{proof}

Denote by $D^\g_{z,z'}$ the selfadjoint operator that is the closure of the
operator $\mathcal{D}^\g_{z,z'}$.

\begin{Proposition}\label{4.2} {\rm $D^\g_{z,z'}$} has pure continuous
spectrum.
\end{Proposition}

\begin{proof} Note that for $\xi=1$, the relation (\ref{2.G}) turns into
$$
\mathcal{D}^\g_{z,z'}=S_{z,z'}\left(\E-\F-\HH\right).
$$
Observe now that $\E-\F-\HH$ is a nilpotent element in
$i\frak{su}(1,1)\subset\frak{sl}(2,\CC)$. For any $r>0$ the matrix
$r(\E-\F-\HH)$ is conjugated to  $\E-\F-\HH$ by means of an element from
$SU(1,1)$. Therefore, $rD^\g_{z,z'}$ is unitarily equivalent to $D^\g_{z,z'}$.
This implies the claim.
\end{proof}

Actually, one can say more: the spectrum is simple and it fills the whole real
axis. This can be deduced, for instance, from a suitable realization of the
representation $S_{z,z'}$. The eigenfunctions of the operator  $D^\g_{z,z'}$
constituting the continuous spectrum can be written down explicitly, see
\cite{BO-Markov}; they are expressed through the classical Whittaker function.

\begin{Proposition}\label{4.3} As $\xi\to1$, the selfadjoint operators  $D_\zz$
converge to the selfadjoint operator  {\rm $D^\g_{z,z'}$} in the strong
resolvent sense.
\end{Proposition}

The proof is the same as in Proposition \ref{3.3}. Just as in Corollary
\ref{3.4} we deduce from this

\begin{Corollary}\label{4.4} As $\xi\to1$, the projection operators $\Projj_+(D_\zz)$
strongly converge to the projection operator {\rm $\Projj_+(D^\g_{z,z'})$}.
\end{Corollary}

The kernel of the limit projection operator is explicitly described in
\cite{BO-Gamma}, \cite{BO-Markov}. It can be written in the integrable form
(\ref{0.A}), where the functions $A$ and $B$ are expressed through the Gamma
function, which explains the origin of the term ``Gamma kernel''.

\end{document}